\documentclass[reqno,12pt,a4paper]{amsart}

\usepackage{mathrsfs}
\usepackage{amssymb}
\usepackage{amsfonts}
\usepackage{latexsym}
\usepackage{amsthm}

\usepackage{graphicx}
\def\lb{\label}

\newcommand{\er}[1]{\textrm{(\ref{#1})}}

\begin{document}


\renewcommand{\theequation}{\arabic{section}.\arabic{equation}}
\theoremstyle{plain}
\newtheorem{theorem}{\bf Theorem}[section]
\newtheorem{lemma}[theorem]{\bf Lemma}
\newtheorem{corollary}[theorem]{\bf Corollary}
\newtheorem{proposition}[theorem]{\bf Proposition}
\newtheorem{definition}[theorem]{\bf Definition}
\newtheorem{remark}[theorem]{\it Remark}

\def\a{\alpha}  \def\cA{{\mathcal A}}     \def\bA{{\bf A}}  \def\mA{{\mathscr A}}
\def\b{\beta}   \def\cB{{\mathcal B}}     \def\bB{{\bf B}}  \def\mB{{\mathscr B}}
\def\g{\gamma}  \def\cC{{\mathcal C}}     \def\bC{{\bf C}}  \def\mC{{\mathscr C}}
\def\G{\Gamma}  \def\cD{{\mathcal D}}     \def\bD{{\bf D}}  \def\mD{{\mathscr D}}
\def\d{\delta}  \def\cE{{\mathcal E}}     \def\bE{{\bf E}}  \def\mE{{\mathscr E}}
\def\D{\Delta}  \def\cF{{\mathcal F}}     \def\bF{{\bf F}}  \def\mF{{\mathscr F}}
\def\c{\chi}    \def\cG{{\mathcal G}}     \def\bG{{\bf G}}  \def\mG{{\mathscr G}}
\def\z{\zeta}   \def\cH{{\mathcal H}}     \def\bH{{\bf H}}  \def\mH{{\mathscr H}}
\def\e{\eta}    \def\cI{{\mathcal I}}     \def\bI{{\bf I}}  \def\mI{{\mathscr I}}
\def\p{\psi}    \def\cJ{{\mathcal J}}     \def\bJ{{\bf J}}  \def\mJ{{\mathscr J}}
\def\vT{\Theta} \def\cK{{\mathcal K}}     \def\bK{{\bf K}}  \def\mK{{\mathscr K}}
\def\k{\kappa}  \def\cL{{\mathcal L}}     \def\bL{{\bf L}}  \def\mL{{\mathscr L}}
\def\l{\lambda} \def\cM{{\mathcal M}}     \def\bM{{\bf M}}  \def\mM{{\mathscr M}}
\def\L{\Lambda} \def\cN{{\mathcal N}}     \def\bN{{\bf N}}  \def\mN{{\mathscr N}}
\def\m{\mu}     \def\cO{{\mathcal O}}     \def\bO{{\bf O}}  \def\mO{{\mathscr O}}
\def\n{\nu}     \def\cP{{\mathcal P}}     \def\bP{{\bf P}}  \def\mP{{\mathscr P}}
\def\r{\rho}    \def\cQ{{\mathcal Q}}     \def\bQ{{\bf Q}}  \def\mQ{{\mathscr Q}}
\def\s{\sigma}  \def\cR{{\mathcal R}}     \def\bR{{\bf R}}  \def\mR{{\mathscr R}}
                \def\cS{{\mathcal S}}     \def\bS{{\bf S}}  \def\mS{{\mathscr S}}
\def\t{\tau}    \def\cT{{\mathcal T}}     \def\bT{{\bf T}}  \def\mT{{\mathscr T}}
\def\f{\phi}    \def\cU{{\mathcal U}}     \def\bU{{\bf U}}  \def\mU{{\mathscr U}}
\def\F{\Phi}    \def\cV{{\mathcal V}}     \def\bV{{\bf V}}  \def\mV{{\mathscr V}}
\def\P{\Psi}    \def\cW{{\mathcal W}}     \def\bW{{\bf W}}  \def\mW{{\mathscr W}}
\def\o{\omega}  \def\cX{{\mathcal X}}     \def\bX{{\bf X}}  \def\mX{{\mathscr X}}
\def\x{\xi}     \def\cY{{\mathcal Y}}     \def\bY{{\bf Y}}  \def\mY{{\mathscr Y}}
\def\X{\Xi}     \def\cZ{{\mathcal Z}}     \def\bZ{{\bf Z}}  \def\mZ{{\mathscr Z}}
\def\O{\Omega}

\newcommand{\gA}{\mathfrak{A}}
\newcommand{\gB}{\mathfrak{B}}
\newcommand{\gC}{\mathfrak{C}}
\newcommand{\gD}{\mathfrak{D}}
\newcommand{\gE}{\mathfrak{E}}
\newcommand{\gF}{\mathfrak{F}}
\newcommand{\gG}{\mathfrak{G}}
\newcommand{\gH}{\mathfrak{H}}
\newcommand{\gI}{\mathfrak{I}}
\newcommand{\gJ}{\mathfrak{J}}
\newcommand{\gK}{\mathfrak{K}}
\newcommand{\gL}{\mathfrak{L}}
\newcommand{\gM}{\mathfrak{M}}
\newcommand{\gN}{\mathfrak{N}}
\newcommand{\gO}{\mathfrak{O}}
\newcommand{\gP}{\mathfrak{P}}
\newcommand{\gQ}{\mathfrak{Q}}
\newcommand{\gR}{\mathfrak{R}}
\newcommand{\gS}{\mathfrak{S}}
\newcommand{\gT}{\mathfrak{T}}
\newcommand{\gU}{\mathfrak{U}}
\newcommand{\gV}{\mathfrak{V}}
\newcommand{\gW}{\mathfrak{W}}
\newcommand{\gX}{\mathfrak{X}}
\newcommand{\gY}{\mathfrak{Y}}
\newcommand{\gZ}{\mathfrak{Z}}

\def\ve{\varepsilon} \def\vt{\vartheta} \def\vp{\varphi}  \def\vk{\varkappa}

\def\Z{{\mathbb Z}} \def\R{{\mathbb R}} \def\C{{\mathbb C}}  \def\K{{\mathbb K}}
\def\T{{\mathbb T}} \def\N{{\mathbb N}} \def\dD{{\mathbb D}} \def\dS{{\mathbb S}}
\def\B{{\mathbb B}}


\def\la{\leftarrow}              \def\ra{\rightarrow}     \def\Ra{\Rightarrow}
\def\ua{\uparrow}                \def\da{\downarrow}
\def\lra{\leftrightarrow}        \def\Lra{\Leftrightarrow}


\def\lt{\biggl}                  \def\rt{\biggr}
\def\ol{\overline}               \def\wt{\widetilde}
\def\no{\noindent}


\newcommand{\fr}{\frac}
\newcommand{\tf}{\tfrac}

\let\ge\geqslant                 \let\le\leqslant
\def\lan{\langle}                \def\ran{\rangle}
\def\/{\over}                    \def\iy{\infty}
\def\sm{\setminus}               \def\es{\emptyset}
\def\ss{\subset}                 \def\ts{\times}
\def\pa{\partial}                \def\os{\oplus}
\def\om{\ominus}                 \def\ev{\equiv}
\def\iint{\int\!\!\!\int}        \def\iintt{\mathop{\int\!\!\int\!\!\dots\!\!\int}\limits}
\def\el2{\ell^{\,2}}             \def\1{1\!\!1}
\def\sh{\sharp}
\def\wh{\widehat}
\def\bs{\backslash}
\def\na{\nabla}

\def\sh{\mathop{\mathrm{sh}}\nolimits}
\def\all{\mathop{\mathrm{all}}\nolimits}
\def\Area{\mathop{\mathrm{Area}}\nolimits}
\def\arg{\mathop{\mathrm{arg}}\nolimits}
\def\const{\mathop{\mathrm{const}}\nolimits}
\def\det{\mathop{\mathrm{det}}\nolimits}
\def\diag{\mathop{\mathrm{diag}}\nolimits}
\def\diam{\mathop{\mathrm{diam}}\nolimits}
\def\dim{\mathop{\mathrm{dim}}\nolimits}
\def\dist{\mathop{\mathrm{dist}}\nolimits}
\def\Im{\mathop{\mathrm{Im}}\nolimits}
\def\Iso{\mathop{\mathrm{Iso}}\nolimits}
\def\Ker{\mathop{\mathrm{Ker}}\nolimits}
\def\Lip{\mathop{\mathrm{Lip}}\nolimits}
\def\rank{\mathop{\mathrm{rank}}\limits}
\def\Ran{\mathop{\mathrm{Ran}}\nolimits}
\def\Re{\mathop{\mathrm{Re}}\nolimits}
\def\Res{\mathop{\mathrm{Res}}\nolimits}
\def\res{\mathop{\mathrm{res}}\limits}
\def\sign{\mathop{\mathrm{sign}}\nolimits}
\def\span{\mathop{\mathrm{span}}\nolimits}
\def\supp{\mathop{\mathrm{supp}}\nolimits}
\def\Tr{\mathop{\mathrm{Tr}}\nolimits}
\def\BBox{\hspace{1mm}\vrule height6pt width5.5pt depth0pt \hspace{6pt}}
\def\where{\mathop{\mathrm{where}}\nolimits}
\def\as{\mathop{\mathrm{as}}\nolimits}


\newcommand\nh[2]{\widehat{#1}\vphantom{#1}^{(#2)}}
\def\dia{\diamond}

\def\Oplus{\bigoplus\nolimits}



\def\qqq{\qquad}
\def\qq{\quad}
\let\ge\geqslant
\let\le\leqslant
\let\geq\geqslant
\let\leq\leqslant
\newcommand{\ca}{\begin{cases}}
\newcommand{\ac}{\end{cases}}
\newcommand{\ma}{\begin{pmatrix}}
\newcommand{\am}{\end{pmatrix}}
\renewcommand{\[}{\begin{equation}}
\renewcommand{\]}{\end{equation}}
\def\eq{\begin{equation}}
\def\qe{\end{equation}}
\def\[{\begin{equation}}
\def\bu{\bullet}

\title[{Trace formulas for  Schr\"odinger operators with complex
potentials}] {Trace formulas for Schr\"odinger operators with
complex potentials}

\date{\today}

\author[Evgeny Korotyaev]{Evgeny Korotyaev}
\address{Saint-Petersburg State University, Universitetskaya nab. 7/9,
St. Petersburg, 199034, Russia, \ korotyaev@gmail.com, \
e.korotyaev@spbu.ru }


\begin{abstract}
\no We consider 3-dim Schr\"odinger operators with a complex
potential. We obtain new trace formulas. In order to prove these
results we study analytic properties of a modified Fredholm
determinant. In fact we reformulate spectral theory problems  as the
problems of analytic functions from Hardy spaces in upper
half-plane.
\end{abstract}

\maketitle






\section {Introduction and main results}
\setcounter{equation}{0}

\subsection{ Introduction.} We consider the Schr\"odinger operator
$ H=-\D+V$ on $L^2(\R^3)$, where the potential $V$ is complex and
satisfies
\[
\lb{1.1}
\begin{aligned}
V\in L^{3\/2}(\R^3)\cap L^2(\R^3).
\end{aligned}
\]
Here $L^p(\R^d), p\ge 1$ is the space equipped with the norm
$\|f\|_p^p=\int_{\R^d}|f(x)|^pdx$ and let $\|\cdot\|=\|\cdot\|_2$.
It is well-known that the operator $H$ has essential  spectrum
$[0,\iy)$ plus $N\in [0, \iy]$ eigenvalues (counted with
multiplicity) in the cut domain $\C\sm [0,\iy)$.  Denote by $\l_j\in
\C\sm [0,\iy), j=1,...,N$ the eigenvalues (according multiplicity)
of the operator $H$. Note, that the multiplicity of each eigenvalue
$\ge 1$, but we call the multiplicity of the eigenvalue its
algebraic multiplicity.

Define  an operator-valued functions $Y_0(k)$ by
\[
\lb{Y0} Y_0(k)=V_1 R_0(k)V_2,   \qqq k\in\C_+,
\]
where  $R_0(k)=(-\D-k^2)^{-1}$ is the free resolvent and we have
used the factorization
\[
\lb{fV}
 V=V_1V_2,  \qqq \where \ \qqq V_1= |V|^{1\/2},\qqq  V_2=|V|^{1\/2}e^{i\arg V}.
\]

 Let $\cB$ denote the class of bounded operators. Let $\cB_1$
and $\cB_2$ be the trace and the Hilbert-Schmidt class equipped with
the norm $\|\cdot \|_{\cB_1}$ and $ \|\cdot \|_{\cB_2}$
correspondingly.  It is well known that $Y_0(k)\in \cB_2$   but
$Y_0(k)\notin \cB_1$ for all $k\in \C_+$. In this case we can not
directly define the determinant $\det (I+Y_0(k))$ and we need some
modification.  It is well-known that the mapping $Y_0: \C_+\to
\cB_2$ is analytic and continuous up to the real line. Then we can
define the modified determinant
\[
\p(k)=\det \lt[(I+ Y_0(k))e^{- Y_0(k)}  \lt], \ \ \ k\in \C_+.
\]
The function $\p(k)$ is analytic in $\C_+$ and continuous up to the
real line. It has $N\le \iy$ zeros (counted with multiplicity)
$k_1,..,k_{N}\in \C_+$ given by $k_n=\sqrt{\l_n}\in \C_+$ and
labeled by
\[
\lb{Lkn} \Im k_1\ge \Im k_2\ge \Im k_3\ge... \ge \Im k_n\ge ...
\]
and define the set $\s_d=\{k_1,...,k_{N}\in \C_+\}$. There are no
other roots in $\C_+$.
 We describe the basic properties
of the determinant $\p$ as an analytic function in $\C_+$.

\begin{theorem}\lb{T1}
Let a potential $V\in L^{3\/2}(\R^3)\cap L^2(\R^3)$. Then the
determinant $\p$ and the function $\p_2={1\/2}\Tr Y_0^2(k)$ are
analytic in $\C_+$ and  continuous up to the real line and
satisfies:
\begin{equation}
\label{D1}
\begin{aligned}
|\p(k)|\le \exp \Big(C_\bu\|V\|_{3\/2}^2\Big),\qqq \forall \ k\in
\ol\C_+,
\end{aligned}
\end{equation}
 where $C_\bu={1\/ 8(4\pi)^{2\/3}}$,
\[
\label{asD}
\begin{aligned}
& \p_2(k)=1+o(1),
\\
& \log \p(k)-\p_2(k)=O(k^{-{1\/2}}) \qqq
\end{aligned}
\]
as $ |k| \to \infty$, uniformly with respect to ${\arg}\,k \in
[0,\pi]$,
\[
\label{epe}
\begin{aligned}
& |\p_2(k)|\le C_\bu\|V\|_{3\/2}^2, \qqq \forall \ k\in \ol\C_+,
\\
& \|\p_2(\cdot +i0)\|\le  {1\/8(4\pi)^{1\/3}}\|V\|_{3\/2}\|V\|,
\end{aligned}
\]
and the function $h(k)=\log |\p(k+i0)|, k\in \R$ satisfies
\[
\label{em1} h\in L^1(-r,r), \qqq h\in L^\a(\R\sm (-r,r))
\]
for some $r>0, \a\in (2,\iy)$ and the zeros $\{k_j\}$ of $\p$ in
$\C_+$   satisfy
\[
\lb{B1} \sum _{j=1}^N \Im k_j<\iy.
\]
\end{theorem}

In fact, $\p$ is a function from the Hardy class. We now introduce
some notations to make this precise.  Let a  function $F(k),
k=u+iv\in \C_+$ be analytic on $\C_+$. For $0<p\le \iy$ we say $F$
belongs the Hardy space $ \mH_p=\mH_p(C_+)$ if $F$ satisfies
$\|F\|_{\mH_p}<\iy$, where $\|F\|_{\mH_p}$ is given by
$$
\|F\|_{\mH_p}=\ca
\sup_{v>0}{1\/2\pi}\rt(\int_\R|F(u+iv))|^pdu\rt)^{1\/p} &
if  \qqq 0< p<\iy\\
 \sup_{k\in \C_+}|F(k)| & if \qqq p=\iy\ac .
$$
Note that the definition of the Hardy space $\mH_p$ involves all
$v=\Im k>0$.

Thus we have that $\p, \p_2\in \mH_\iy$ and $\p_2\in \mH_2$.

\subsection{Trace formulas and estimates}
Theorem \ref{T1} shows the basic analytic properties of the function
$\p$, in particular, $\p\in \mH_\iy$. In order to study zeros of
$\p(k)$ in the upper-half plane we need to study the Blaschke
product, defined by
\begin{equation}
\label{Bk}
 B(k)=\prod_{j=1}^N\rt(\frac{k-k_j}{ k-\ol k_j}\rt),\qqq k\in \C_+.
\end{equation}
where the Blaschke product $B(k)$ converges absolutely for $k\in
\C_+$, since \er{B1} holds true. Moreover, we have the function
$B\in \mH_\iy$ with $\|B\|_{\mH_\iy}\le 1$,  (see e.g. \cite{Ko98}
or \cite{G81}).

We describe the canonical factorization of $\p$.

\begin{theorem}
\lb{T2} Let $V\in L^{3\/2}(\R^3)\cap L^2(\R^3)$. Then $\p$ has a
canonical factorization in $\C_+$ given by
\[
\lb{facp}
\begin{aligned}
\p=\p_{in}\p_{out},\qqq &  \p_{in}(k)=B(k)e^{-iK(k)},\qqq  K(k)=
{1\/\pi}\int_\R {d\n(t)\/k-t},
\\
& \p_{out}(k)=e^{iM(k)},\qqq  M(k)= {1\/\pi}\int_\R
 {\log |\p(t)|\/k-t}dt.
 \end{aligned}
\]
\no $\bu$  $d\n(t)\ge 0$ is a singular compactly supported measure
on $\R$ and for some $r_*>0$ satisfies
\[
\lb{smsx}
\begin{aligned}
\int_\R d\n(t)<\iy,\qqq
 \supp \n \ss \{z\in \R: \p(z)=0\}\ss [-r_*, r_*].
 \end{aligned}
\]

\no $\bu$  The function $K$ has an analytic continuation from $\C_+$
into the domain $\C\sm [-r_*, r_*]$   and has the following Taylor
series
\[
\lb{} K(k)=\sum_{j=0}^\iy {K_j\/k^{j+1}},\qqq K_j={1\/\pi}\int_\R
t^jd\n(t).
\]
\no $\bu$ The Blaschke product $B$ has an analytic continuation from
$\C_+$ into the domain $\{|k|>r_0\}$, where $r_0=\sup |k_n|$ and has
the following Taylor series
\[
\begin{aligned}
\lb{B6}
\log  B(z)=-i{B_0\/k}-i{B_1\/2k^2}-i{B_2\/3k^3}-... \qqq as \qqq |k|>r_0,\\
 B_0=2\sum_{j=1}^N\Im k_j,\qqq B_1=2\sum_{j=1}^N\Im k_j^2,..., \
B_n=2\sum_{j=1}^N\Im k_j^{n+1},.....
\end{aligned}
\]
where each sum $B_n, n\ge 1$ is absolutely convergence and satisfies
\[
|B_n|\le 2\sum_{j=1}^N|\Im k_j^n|\le {\pi\/2}(n+1)r_0^{n}B_0.
\]

\end{theorem}

{\bf Remarks.} 1) The function $\p_{in}$ is the inner factor of $\p$
and the function $\p_{out}$ is the outer factor of $\p$. These
results will be used in the proof of trace formulas.

2) Due to \er{em1} the integral $M(k),k\in \C_+$ in \er{facp}
converges absolutely.

\begin{corollary}
\lb{T2x} Let $V\in L^{3\/2}(\R^3)\cap L^2(\R^3)$. Then  the
following trace formula
\[
\lb{tre1} -2k\Tr \rt(R(k)-R_0(k) +R_0(k)VR_0(k) \rt)= \sum {2i\Im
k_j\/(k-k_j)(k-\ol k_j)}+{i\/\pi}\int_{\R}{d\m(t)\/(t-k)^2},
\]
holds true for any $k\in \C_+\sm\s_d$, where  the measure
$d\m(t)=\log |\p(t)|dt-d\n(t)$ and the series converges uniformly in
every bounded disc in $\C_+\sm \s_d$. If in addition $V\in
L^1(\R^3)$, then
\[
\lb{tre2} -2k\Tr \rt(R(k)-R_0(k)\rt)  +{i\/4\pi}\int_{\R^3}V(x)dx=
\sum {2i\Im k_j\/(k-k_j)(k-\ol
k_j)}+{i\/\pi}\int_{\R}{d\m(t)\/(t-k)^2}.
\]

\end{corollary}

In order to show trace formulas we introduce the space $W_m$  by
\[
W_m=\lt\{f\in L^2(\R^3):  |\pa_x^\a f(x)|\le
{C_V\/(1+|x|)^{\b+|\a|}}, \ \forall \ |\a|\le 2m+1 \rt\},\qq m\geq
1,
\]
where $\b>3$ and $m\ge 1$ is an integer. If $V\in W_{m}$, then the
function $\p(\cdot)$ is analytic in $\C_+$ and continuous up to the
real line and  satisfies
\[
\lb{apm}
  \log \p(k)={Q_0\/ik }+{Q_2\/ik ^{3}}+{Q_4\/ik^{5}}+\dots
  +{Q_{2m}+o(1)\/ik^{2m+1}},
\]
 as $|k |\to \iy$ uniformly in $\arg k \in [0,\pi]$ (see \cite{B66},
and also \cite{C81,G85,P82,R91}), where
\[
\lb{Qn}
Q_0={1\/ 16\pi }\int_{\R^3} V^2(x)dx,  \ \ \ \ \
Q_2={1\/ 3\pi 4^3}\int_{\R^3} \rt( (\nabla V(x))^2+2V^3(x)\rt)dx,.....
\]
 In Theorem \ref{T3} we show that  the function
$M(k), k\in \C_+$ defined by \er{facp} satisfies
\[
\begin{aligned}
\lb{Kasm} M(k)={1\/\pi}\int_\R {
h(t)\/k-t}dt={\cJ_0-iI_{0}\/t}+{\cJ_1\/t^2}+...
+{\cJ_{2m}-iI_{2m}+o(1)\/t^{2m+1}},\\
\end{aligned}
\]
where $h:=h_{-1}:=\log |\p(\cdot)|,$ as $\Im k\to \iy$ and
\[
\begin{aligned}
\lb{Kas}
 I_j=\Im Q_j, \qq h_{j}=t^{j+1}(h-P_{j}) ,\qq
P_j={I_{0}\/t}+{I_{2}\/t^3}+...+{I_{j}\/t^{j+1}}, \qq
\\
\cJ_0={\rm v.p.}{1\/\pi}\int_\R h(t) dt,\qqq \cJ_1={\rm
v.p.}{1\/\pi}\int_\R \big(th(t)-I_0\big)dt,\qq \cJ_j={\rm
v.p.}{1\/\pi}\int_\R h_{j-1}(t)dt,
\end{aligned}
\]
$j=0,1,...,2m$. Here $I_{2j+1}=0$ and all integrals in \er{Kas}
converges, since $\p$ satisfies \er{apm}.

\begin{theorem}\label{T3}
({\bf Trace formulas}) Let $V\in L^{3\/2}(\R^3)\cap L^2(\R^3)$. Then
\[
\begin{aligned}
\lb{tr12} B_0+{\n(\R)\/\pi }={1\/\pi}\int_\R\log |D_4(k)|dk.
\end{aligned}
\]
 Let in addition a potential $V\in W_{m+1}$ for some $m\ge 1$. Then
the function $M$ defined by \er{facp} satisfies \er{Kasm}. Moreover,
the following identities hold true:
\[
\lb{trj}
\begin{aligned}
& {B_1\/2}+K_1=\cJ_1,\\
& {B_j\/j+1}+K_j=\Re Q_j+\cJ_j,\qq j=2,3,...,2m.
\end{aligned}
\]
\end{theorem}

\medskip
\noindent {\bf Remark.} 1) There are a lot of results about the
trace formulas.  The trace formulas similar to (\ref{tr12}),
(\ref{trj})  were proved by Buslaev \cite{B66} for real potentials,
see also \cite{C81} and \cite{G85,P82,R91}. Trace formulas for the
case Stark operators and magnetic Schr\"odinger operators are
discussed in \cite{KP03}, \cite{KP04}. Trace formulas for
Schr\"odinger operators on the lattice are considered in \cite{IK12}
and for the case of complex potentials in
 \cite{K17}, \cite{KL16}. In our consideration we use the methods
from \cite{K17}, \cite{KL16}.

\no 2) In the case of  complex potentials there is an additional
term associated with the singular measure $\n$, see \er{tr12},
\er{trj}.

\noindent 3) Higher regularity of  $V$ implies more trace formulas.

\noindent 4) Buslaev mainly considered the phase $\phi_{sc}=\arg
\p$. In the present paper, the trace formulas for the conjugate
function $\log|\p(k)|$ are proved.

  \begin{theorem} ({\bf Estimates})
\lb{T4} Let $V\in L^{3\/2}(\R^3)\cap L^2(\R^3)$. Then the following
estimate
\[
\begin{aligned}
\lb{eBs}
 & {\n(\R)}+\sum_{j=1}^N\Im k_j\le \|V\|^2F(\|V\|_{3\/2}),
\end{aligned}
\]
holds true, where
\[
\begin{aligned}
\lb{eBsf} & F(\l)=a_1 \l^{1\/2}+a_2 \l +a_3 \l^{3\/2},\l>0,
\\
& a_1=C_2{(68)^{1\/4}\/(4\pi)^{7\/6}},\qq a_2={\sqrt
{68}\/(4\pi)^{5\/6}}C_2^2,\qqq a_3=C_2^3{(68)^{1\/4}\/ 3\pi
(4\pi)^{1\/2}}
\end{aligned}
\]
and  $C_2$ is the constant from \er{F2}.
  \end{theorem}

Recently uniform bounds on eigenvalues of Schr\"odinger operators in
$\R^d$ with complex-valued potentials decaying at infinity attracted
attention of many specialists. For example, bounds on single
eigenvalues were proved  in \cite{DN02,F11,Sa10} and bounds on sums
of powers of eigenvalues were found in
\cite{FLLS06,LS09,DHK09,DHK13, FS17,F15}. These bounds generalise
the Lieb--Thirring bounds \cite{LT76} to the non-selfadjoint
setting. Note that in \cite{FS17} (Theorem 16) the authors obtained
estimates on the sum of the distances between the complex
eigenvalues and the continuous spectrum $[0,\infty)$ in terms of
$L^p$-norms of the potentials. Note that almost no results are known
on the number of eigenvalues of Schr\"odinger operators with complex
potentials. We referee here to a recent paper \cite{FLS16} where the
authors discussed this problem in details in odd dimensions. In
\cite{Sa10} Safronov obtained some trace formulae in the continuous
case for showing that the series of imaginary parts of square roots
of eigenvalues is convergent.

Finally we note Schr\"odinger operators with complex potentials on
the lattice was studied in \cite{K17},  \cite{KL16}. Here trace
formulas and different estimates of complex eigenvalues in terms of
potentials were determined and the Hardy space in the unit disk was
used. In the present paper for the continuous case we use the
technique from \cite{K17}, \cite{KL16}, but it is more natural to
use the Hardy space in the upper half-plane.


\section {Preliminaries }
\setcounter{equation}{0}

\subsection{Determinants and Trace class operators}
Let us recall some well-known facts.

\no $\bu$ Let $A, B\in \cB$ and $AB, BA\in \cB_1$. Then
\begin{equation}
\label{AB} {\rm Tr}\, AB={\rm Tr}\, BA,
\end{equation}
\begin{equation}
\label{1+AB} \det (I+ AB)=\det (I+BA).
\end{equation}

\no $\bu$ If $A,B\in \cB_1$, then
\begin{equation}
\label{DA1} |\det (I+ A)|\le e^{\|A\|_{\cB_1}}.
\end{equation}

\no $\bu$ We introduce the space $\cB_p, p\ge 1$  of compact
operators $A$ equipped with the norm
$$
\|A\|^p_{\cB_p}=\Tr (A^*A)^{p\/2}<\infty.
$$

In the case $A\in \cB_n, n\ge 2$ we have that the operator  $(I+
A)e^{-A+\G_n(A)}-I\in \cB_1$, where
 $\G_n(z)=\sum_{2}^{n-1}{1\/j}(-z)^j$.
Thus  we define the  modified determinant $\det_n(I+ A)$ by
\[
\begin{aligned}
\label{DA2}
& \det_2 (I+ A) =\det \rt((I+ A)e^{-A}\rt),\\
& \det_n(I+ A) =\det \rt((I+ A)e^{-A+\G_n(A)} \rt), \qq n\ge 3.
 \end{aligned}
 \]
 Note the  modified  determinant satisfies
\[
\label{DA2xi}
\begin{aligned}
\det_n(I+ A) =e^{\Tr \G_n(A)}\det_2(I+ A),\qqq \ if \qq A\in \cB_2
  \end{aligned}
\]
and  $I+ A$ is invertible if and only if $\det_n(I+ A)\ne 0$, see
Chapter IV in \cite{GK69}).

\no $\bu$ We need the estimate from \cite{KS17}
\[
\label{DA2x}
\begin{aligned}
 |\det_n(I+ A)|\le e^{{c_n\/2}\|A\|_{\cB_n}^n},\qq c_n=\ca 1, &
n=2\\
17 & n\ge 3 \ac .
  \end{aligned}
\]

\no $\bu$  Suppose a function $A(\cdot): \mD  \to \cB_1$ is analytic
for a domain $\mD \subset {\C}$, and the operator $(I+A(z))^{-1} $
is bounded  for any $z\in \mD$. Then for $F(z)=\det (I+A(z))$ we
have
$$
 F'(z)= F(z)\Tr \big(I+A(z)\big)^{-1}A'(z)\qqq \forall \ z\in \mD.
$$

\subsection{Analytic functions}
 The kernel of the free resolvent $R_0(k)=(-\D-k^2)^{-1}$ has the
form
\[
\lb{dR0} R_0(x-y,k)={e^{ik|x-y|}\/ 4\pi |x-y|},\ \ \ \ \ x,y\in
\R^3, \ \ \ \ k\in \C_+.
\]
Recall that  $Y_0(k)$ is defined in \er{Y0} and  $ \s_d
=\{k_j=\sqrt{\l_j}, j=1,...,N\}\ss \C_+$. Define an operator
$Y(k)=|V|^{1/2} R(k)V^{1\/2}$ for $k\in\C_+$. Below we will use the
standard identity
\begin{equation}
 (I+Y_0(k))(I-Y(k))=I,\quad \forall \
 k\in {\C}_+\sm \s_d.
\label{S2I+Q0I+Q}
\end{equation}


Below we need the Hardy-Littlewood-Sobolev inequality for $f,g\in
L^{3\/2}(\R^3)$:
\[
\lb{HLS1} \int_{\R^{6}}{|f(x)||g(y)|\/|x-y|^2}dxdy\le
\pi^{4\/3}4^{1\/3}\|f\|_{3\/2}\|g\|_{3\/2}.
\]
We recall the well known results.

\begin{lemma}
\lb{T21} i) Let $V\in L^{3\/2}(\R^3)$. Then the function
$Y_0:\C_+\to \cB_2$ is analytic and is continuous up to the real
line and satisfies for all $k\in\ol\C_+$:
\begin{equation}
\lb{B22} \|Y_0(k)\|_{\cB_2}^2\le \int_{R^6}{|V(x)||V(y)|\/
(4\pi)^2|x-y|^2}dxdy \le 2C_\bu \|V\|_{3\/2}^2,\qqq C_\bu=
{1\/8(4\pi)^{2\/3}}.
\end{equation}
ii) If $V\in L^2(\R^3)$ and $\Im k>0$, then
\begin{equation}
\lb{B21}
\begin{aligned}
& \|Y_0(k)\|_{\cB_2}\le {\|V\|_2\/ \sqrt{8\pi \Im k}},
\\
&   VR_0(k)\in \cB_2.
\end{aligned}
\end{equation}
\no iii) If $V\in L^1(\R^3)$, then the function $Y_0':\C_+\to \cB_1$
is analytic and is continuous up to the real line and satisfies for
all
 $k\in {\C}_+$:
\begin{equation}
\lb{B24} {\rm Tr}\ Y_0'(k)= 2k{\rm
Tr}\,\big(VR_0^2(k)\big)={i\/4\pi}\int_{\R^3}V(x)dx.
\end{equation}

\end{lemma}

{\bf Proof}. i) The Hardy-Littlewood-Sobolev inequality \er{HLS1}
and \er{dR0} yield
$$
\|Y_0(k)\|_{\cB_2}^2={1\/
(4\pi)^2}\int_{R^6}|V(x)||V(y)|{e^{-2v|x-y|}\/ |x-y|^2}dxdy\le {1\/
(4\pi)^2}\int_{R^6}{|V(x)||V(y)|\/ |x-y|^2}dxdy
$$
$$
\le {\pi^{4\/3}4^{1\/3}\/
(4\pi)^2}\|V\|_{3\/2}^2=2C_\bu\|V\|_{3\/2}^2.
$$
Similar arguments give well-known results that the operator-valued
functions
 $Y_0:\C_+\to \cB_2$ is analytic and is continuous up to the real
line.

ii) Using  the identity \er{dR0}, $v=\Im k>0$ and the Schwartz
inequality we obtain
$$
\|Y_0(k)\|_{\cB_2}^2={1\/
(4\pi)^2}\int_{R^6}|V(x)||V(y)|{e^{-2v|x-y|}\/ |x-y|^2}dxdy\le
%
\frac{1}{(4\pi)^2}\int_{R^6}|V(x)|^2{e^{-2v|x-y|}\/ |x-y|^2}dxdy
$$
$$
={\|V\|_2^2\/ (4\pi)^2} \int_{R^3}{e^{-2v|y|}\/|y|^2}dy
={\|V\|_2^2\/ (4\pi)}  \int_{0}^\iy{e^{-2vr}}dr={\|V\|_2^2\/ 8\pi
v}.
$$

iii) It  $V\in L^1(\R^3)$, then similar arguments give that the
function $Y_0'(k)$ is analytic in $\C_+$ and is continuous up to the
real line.   Using \er{dR0}, we obtain
$$
{\rm Tr}\, Y_0'(k)= 2k{\rm Tr}\,(VR^2(k))=\frac{i}{4\pi }\int_{{\
R}^3} V(x)dx,\qqq  k\in {\C}_+.
$$
 \BBox

Below we need  the following estimates from  \cite{FS17}:

{\it Let  ${3\/2}<q\le 2$ and $p={2q\/3-q}$.   Then the
${\cB}_p$-norms of the operator $Y_0(\lambda)$ satisfies
\[
\lb{F1} \|Y_0(\lambda)\|_{{\cB}_p}\le C_q |k|^{{3\/q}-2}\|V\|_{{q}},
\]
\[
\lb{F2}   \|Y_0(k)\|_{{\cB}_4}\leq C_2 |k|^{-{1\/2}}\|V\|_2,
\]
where the constant $C_q>0$ depends on $q$ only.}

\medskip

If $V\in L^{3\/2}(\R^3)$, then we have $\|Y_0(k)\|$=o(1) as $|k|\to
\iy$ uniform with respect to ${\arg}\,k \in [0,\pi]$.

In this case we can define the radius $r_0>0$ by the condition
\[
\lb{dr0} \sup_{\Im k\ge 0, |k|\ge r_0} \|Y_0(k)\|\le {1\/2}.
\]

Let $\p_2={1\/2}\Tr Y_0^2$. Define the Fourier transformation
$$
V(x)={1\/(2\pi)^{3\/2}}\int_{\R^3}e^{i(p,x)}\hat V(p)dp,
$$
where $(\cdot,\cdot)$ is the scalar product in $\R^3$.

\begin{lemma}
\lb{Tdt1} Let $V\in L^{3\/2}(\R^3)\cap L^2(\R^3)$. Then the function
$\p_2={1\/2}\Tr Y_0^2$ satisfies
\[
\lb{4.2} \p_2(k)=\int_0^\iy e^{i2kt}\g(t)dt=-C\int_{\dS^2}d\o
\int_{\R^3}{\hat V(p)\hat V(-p)\/i(2k-(p,\o))}dp,
\]
where $C={1\/2(4\pi )^2}$ and
\[
\lb{4.2b} \g(t)=C\int_{|\o |=1}d\o \int_{\R^3}V(x-t\o
)V(x)dx=C\int_{\dS^2}d\o \int_{\R^3}\hat V(p)\hat
V(-p)e^{-it(p,\o)}dp,
\]
and
\[
\label{D1b}
\begin{aligned}
 \|\g\|_1\le C_\bu\|V\|_{3\/2}^2 \qqq & if \qq
V\in L^{3\/2}(\R^3),\\
\|\g\|_{\iy}\le { \|V\|^2\/8\pi} \qqq & if \qq V\in L^2(\R^3),
\\
\|\p_2\|\le  {1\/8(4\pi)^{1\/3}}\|V\|_{3\/2}\|V\| \qqq & if \qq V\in
L^2(\R^3)\cap L^{3\/2}(\R^3),
\end{aligned}
\]
\[
\label{D1c} \p_2(k)=-{Q_0+o(1)\/ik} \qqq \as \ \ \Im k\to 0,
\]
\[
\label{D1j} \p_2\in \mH_2.
\]
\[
\label{asp2} \p_2(k)=o(1) \qqq \as \ \ |k|\to \iy,
\]
uniform with respect to ${\arg}\,k \in [0,\pi]$, and here
$Q_0={1\/16\pi}\int_{\R^3}V^2(x)dx$.

\end{lemma}

\no {\bf  Proof.}
 Due to \er{dR0} we have
$$
 \p_2(k)={1\/2}\int_{\R^6}V(x)R_0^2(x-y,k)V(y)dxdy=
C\int_{\R^6}{V(x)V(y)e^{i2k|x-y|}\/  |x-y|^2}dxdy.
$$
 If we set $x-y=t\o$, where $\o$ belongs to the unit sphere $\dS^2$ and $t=|x-y|>0$, then we obtain
\[
\lb{4.27} \p_2(k)=\int_0^\iy e^{i2kt}\g(t)dt, \ \ \ \ \
\g(t)=C\int_{\dS^2}d\o \int_{\R^3}V(x-t\o )V(x)dx.
\]
Using the Fourier transformation  we rewrite the function $\g$ in
the form
\[
\lb{4.28}
\begin{aligned}
\g(t)=C\int_{\dS^2}d\o \int_{\R^3}V(x-t\o )V(x)dx=C\int_{\dS^2}d\o
\int_{\R^3}\hat V(p)\hat V(-p)e^{-it(p,\o)}dp.
\end{aligned}
\]
Thus we have for $k\in \C_+$:
\[
\lb{4.2z} \p_2(k)=C\int_{\dS^2}d\o \int_{\R^3}\hat V(p)\hat V(-p)dp
\int_0^\iy e^{it(2k-(p,\o))}dt=C\int_{\dS^2}d\o \int_{\R^3}{\hat
V(p)\hat V(-p)\/i((p,\o)-2k)}dp.
\]
From this identity and the Lebesgue Theorem we have
$$
\begin{aligned}
\p_2(k)+{Q_0\/ik}=C\int_{\dS^2}d\o \int_{\R^3}\hat V(p)\hat
V(-p)\rt({1\/i2k}-{1\/i(2k-(p,\o))}\rt)dp
\\
={C\/i2k}\int_{\dS^2}d\o \int_{\R^3}\hat V(p)\hat V(-p) {p\o
dp\/2k-p\o}={o(1)\/k},
\end{aligned}
$$
since ${p\o \/2k-p\o}\to 0$ as $\Im k\to \iy$ and $V\in
L^{2}(\R^3)$.

We show \er{D1b}. Let $V\in L^{3\/2}(\R^3)$. Then using \er{B22} we
have
$$
\int_0^\iy |\g(t)|dt\le C\int_0^\iy \int_{\dS^2}dt d\o
\int_{\R^3}|V(x-t\o )V(x)|dx=C\int_{R^6}{|V(x)||V(y)|\/
|x-y|^2}dxdy<C_\bu\|V\|_{3\/2}^2.
$$
If $V\in L^{2}(\R^3)$, then for any $t\ge 0$ we obtain
\[
\lb{4.2x} |\g(t)|\le C\int_{\dS^2} d\o \int_{\R^3}|V(x-t\o
)V(x)|dx\le C|\dS^2|\cdot \|V\|^2={\|V\|^2\/8\pi}.
\]
Thus $\p_2\in L^2(\R)$ and $\p_2\in \mH_2$  and we have
$$
{1\/\pi}\|\p_2\|^2=\|\g\|^2\le\|\g\|_\iy\int_0^\iy |\g(t)|dt\le
{\|V\|^2\/8\pi}C_\bu\|V\|_{3\/2}^2,
$$
which give \er{D1b}. Using the identity \er{4.2}, where $\g\in
L^1(0,\iy)$ we obtain \er{asp2}.
 \BBox

\begin{lemma}
\lb{Tp3} Let $V\in L^{3\/2}(\R^3)\cap L^2(\R^3)$. Then the function
$\p_3={1\/3}\Tr Y_0^3$  is analytic in $\C_+$ and is continuous up
to the real line and satisfies:
\[
\lb{p23}
\begin{aligned}
& |\p_3(k)|\le {1\/96\pi}\|V\|_{3\/2}^3,
\\
& |\p_3(k)|\le {1\/3}\|Y_0^3(k)\|_{\cB_1}\le
{C\/|k|^{1\/2}}\|V\|_{3\/2}^2\|V\|,\qq \forall \
k\in\ol\C_+\sm\{0\},
\end{aligned}
\]
\[
\lb{p3as}
\begin{aligned} \p_3(i\t)=O(\t^{-{3\/2}})\qqq \as \qq \t\to +\iy.
\end{aligned}
\]
\end{lemma}

\no {\bf  Proof.} From Lemma  \ref{T21} we deduce that the function
$\p_3={1\/3}\Tr Y_0^3$  is analytic in $\C_+$ and is continuous up
to the real line. From \er{B22} we get
$$
|\p_3(k)|={1\/3}|\Tr Y_0^3(k)|\le {1\/3}\|Y_0(k)\|_{\cB_2}^3\le
{1\/96\pi}\|V\|_{3\/2}^3\qq \forall \ k\in \ol\C_+.
$$
Let $k\in \ol\C_+$  and $k\ne 0$. Then from \er{B22}, \er{F2} we
obtain
$$
\begin{aligned}
& |\p_3(k)|\le {1\/3}\|Y_0^3(k)\|_{\cB_1}\le
{1\/3}\|Y_0(k)\|_{\cB_2}^2\|Y_0(k)\|_{\cB_4}\le
{C\/|k|^{1\/2}}\|V\|_{3\/2}^2\|V\|
\end{aligned}
$$
for some absolute constant $C$. From \er{B21} we have as $k=i\t,
\t\to \iy$:
$$
\begin{aligned}
|\p_3(k)|\le {1\/3}\|Y_0(k)\|_{\cB_2}^3=O(k^{-{3\/2}}).
\end{aligned}
$$
 \BBox

We describe properties of the determinants $\p=D_2$ and $D_4=\det
\rt((I+ Y_0)e^{-Y_0+{1\/2}Y_0^2-{1\/3}Y_0^3}\rt)$.

\begin{lemma}
\lb{TD} Let $V\in L^{3\/2}(\R^3)\cap L^2(\R^3)$ and let $k\in
\ol\C_+$. Then each $D_m\in \mH_\iy, m\ge 2$ and is continuous up to
the real line and satisfies:
\[
\label{D1y} |\p(k)|\le e^{C_\bu \|V\|_{3\/2}^2},
\]
\[
\label{D1a}  \log | D_m(\cdot+i0)|\in L_{loc}^1(\R),
\]
and if in addition $|k|\ge r_0$, where $r_0$ is given by \er{dr0},
then
\begin{equation}
\label{D2}
\begin{aligned}
\log \p(k) =\Tr \G_m(Y_0(k))+\log D_m(k),
\\
\log D_m(k)= \sum _{n=m}^{\infty}{(-1)^{n+1}\/n} \Tr  Y_0^n(k),
\end{aligned}
\end{equation}
where the series converges absolutely and uniformly, and
\begin{equation}
\label{D3}
 \Big|\log D_m(k)\Big|\le {2\/m}\|Y_0^m(k)\|_{\cB_1},
\end{equation}
\[
\label{D4} \p(k)=1+o(1)
\]
as $|k|\to \iy$ uniform with respect to ${\arg}\,k \in [0,\pi]$.

\end{lemma}

{\bf Proof}. Lemma \ref{T21} gives that the operator-valued function
$Y_0: \C_+\to \cB_2$ is analytic in $\C_+$ and is continuous up to
the real line. Then  $\p(k)$ is analytic in $\C_+$ and and is
continuous up to the real line. The estimate \er{D1y} follows from
\er{B22} and \er{DA2x}. Moreover, \er{D1a} holds true,  since  we
have $D_m\in \mH_\iy$,  see e.g. \cite{Ko98}.

 We denote the series in \er{D2} by $F(k)$. Since
\begin{equation}
 |{\rm Tr}\, Y_0^n(k)|\leq \|Y_0^m(k)\|_{\cB_1}\|Y_0(k)\|^{n-m}\leq
\|Y_0^m(k)\|_{\cB_1} \ve_k^{n-m}, \quad
\ve_k=\|Y_0(k)\|\le\frac{1}{2}, \label{S2TrQ0kn}
\end{equation}
The function $F(k)$  converges absolutely and uniformly, and each
term is analytic in $|k|> r_0$. Then $F(k)$  is analytic in $|k|>
r_0$.  Moreover, differentiating \er{D2} and using
(\ref{S2I+Q0I+Q}),
 we have
$$
F'(k)=-i\sum_{n=2}^{\infty}{\rm Tr}\,(-Y_0(k))^{n-1}Y_0'(k)= i{\rm
Tr}\, Y(k)Y_0'(k),  \quad  |k|>C_0.
$$
Then we have $F(k) = - i\log \p(k)$, since $F(i\tau)=o(1)$ as $\tau
\to \infty$. Using \er{D2} and \er{S2TrQ0kn},  we obtain (\ref{D3}).

We show \er{D4}. From \er{D2} we have $\log \p(k) =\p_2(k)+\log
D_3(k),$ as $ |k|>C_0$. Asymptotics (\ref{D3}) and \er{p23} give
$\log D_3(k)=O(k^{-{1\/2}})$ and using \er{asp2} we obtain \er{D4}.
 \BBox

Define the function $\p_j={1\/j}\Tr  Y_0^2(k), j\ge 2$.
 Using \er{DA2xi} we have
\[
\lb{d1a}
\begin{aligned}
D_4=D_2e^{\p_2- \p_3},
\end{aligned}
\]
which yields
\[
\lb{d2x}
\begin{aligned}
h(k)=\log |D_2(k)|=\log |D_4(k)|-\Re (\p_2(k)- \p_3(k)), \qq k\in
\ol\C_+.
\end{aligned}
\]

\begin{lemma}
\lb{Tdtx} Let $V\in L^{3\/2}(\R^3)\cap L^2(\R^3)$. Then  the
following identity holds true:
\[
\lb{M24} {1\/\pi}\int_{\R}{\log
|\p(t)|\/k-t}dt={1\/\pi}\int_{\R}{\log
|D_4(t)|\/k-t}dt+i\p_2(k)-i\p_3(k),\qqq \forall \ k\in\C_+,
\]
\[
\lb{M41} \log |D_4(t+i0)|\in L^1(\R),
\]
and if $k=i\t, \t\to +\iy$, then
\[
\lb{Mas1}
\begin{aligned}
& {1\/\pi}\int_{\R}{\log |\p(t)|\/k-t}dt={C_M+o(1)\/k} ,
\\
& C_M={1\/\pi}\int_{\R}\log |D_4(k)|dt-Q_0.
\end{aligned}
\]
\end{lemma}

\no {\bf  Proof.}  Let $f(t)=\log |D_4(t+i0)|,  t \in \R$. From
Lemma \ref{TD} we have $f\in L_{loc}^1(\R)$. Moreover, from \er{D3},
\er{F2} we deduce that
\begin{equation}
\label{D3xx}
 \Big|\log D_4(t+i0)\Big|\le \|Y_0^4(t+i0)\|_{\cB_1}=O(t^{-2})
 \quad \as \ t\to \pm \iy,
\end{equation}
which gives \er{M41}. Using \er{d2x} and \er{Ki1} and Lemmas
\ref{Tdt1}, \ref{Tp3} we obtain for all $k\in\C_+$:
$$
\begin{aligned}
{1\/\pi}\int_{\R}{h(t)dt\/k-t}={1\/\pi}\int_{\R}{f(t)-\Re
\p_2(t)+\Re \p_3(t)\/k-t}dt
={1\/\pi}\int_{\R}{f(t)\/k-t}dt+i\p_2(k)-i\p_3(k).
\end{aligned}
$$

We show \er{Mas1}. Let $k=i\t, \t\to \iy$.  The Lebesgue Theorem
and \er{M41} give
\[
\lb{asMi}
\begin{aligned}
{1\/\pi}\int_{\R}{f(t)dt\/k-t}-{1\/\pi}\int_{\R}{f(t)dt\/k}={1\/\pi
k}\int_{\R}{tf(t)dt\/k-t}={o(1)\/ k},
\\
\end{aligned}
\]
since ${t\/i\t-t}\to 0$ as $\t\to \iy$ for each $t\in \R$. Then
substituting  asymptotics \er{asMi},\er{D1c}  and \er{p3as} into
\er{M24} yield \er{Mas1}.
 \BBox


{\bf Proof of Theorem \ref{T1}.} Due to Lemmas \ref{TD} the function
$\p(k)$ is analytic in $\C_+$ and satisfies \er{D1}, since we have
\er{B22}.

From \er{asp2} we have the first asymptotics in \er{asD}. From
\er{D2} we have
$$
\log \p(k)-\p_2(k)=\log D_3(k), \qq k\in \ol\C_+.
$$
and from \er{D3}, \er{p23} we get  the second one  in \er{asD}.

From \er{B22} we have the first estimate in \er{epe} and from
\er{D1b} we have the second one.

The function $\p\in \mH_\iy$, then  \er{Ki2} and asymptotics
\er{asD} yield \er{em1}.

Asymptotics \er{asD} gives that all zeros of $\p$ are uniformly
bounded. Then the zeros of $\p$  in $\C_+$ satisfy \er{B1}, since it
is standard for functions from $\mH_\iy$, see p. in \cite{Ko98}.
\BBox

\section {Proof of main theorems }
\setcounter{equation}{0}

We describe the determinant $\p(k), k\in \C_+$ in terms of a
canonical factorization. We remark that (see page 53 in \cite{G81}),
in general, in the upper half plane the condition \er{B1} is replace
by the lager condition:
\[
\lb{BLy} \sum {\Im z_j\/1+|z_j|^2}<\iy,
\]
and the Blaschke product with zeros $z_j$ has the form
\[
\lb{BL2x} B(z)={(z-i)^m\/(z+i)^m} \prod_{z_j\ne
0}^N{|1+z_j^2|\/1+z_j^2}\rt(\frac{z-z_j}{z-\ol z_j}\rt), \qqq z\in
\C_+.
\]
If all  moduli $|z_n|$ are uniformly bounded,  the estimate \er{BLy}
becomes $\sum \Im z_j<\iy$ and the convergence factors in \er{BL2x}
are not needed, since $\prod_{z_j\ne 0}^N\rt(\frac{z-z_j}{z-\ol
z_j}\rt)$ already converges.

{\bf  Proof of Theorem \ref{T2}.} From Theorem \ref{T1} we deduce
that the determinant $\p\in \mH_\iy$.  Recall that due to \er{D4}
each zero $\in \C_+$ of $ \p$  belongs  to the half-disc $\{k\in
\C_+: |k|\le r_*\}$ for some $r_*>0$. From Theorem \ref{T1} we
deduce that the function $\p\in \mH_\iy^m$ satisfies all conditions
 from Theorem \ref{Tcf}  and  we obtain all results in Theorem \ref{T2}.

Note that all needed results about the Blaschke product $B$ are
proved in Lemma \ref{TY2}. \BBox

We prove the first main result about the trace formulas.

{\bf  Proof of Corollary  \ref{T2x}.} From \er{facp} we obtain for
all $k\in \C_+$:
\[
\begin{aligned}
\lb{derfx}
{\p'(k)\/\p(k)}={B'(k)\/B(k)}-{i\/\pi}\int_\R{d\m(t)\/(k-t)^2},\qq
\qqq {B'(k)\/B(k)}=\sum {2i\Im k_j\/(k-k_j)(k-\ol k_j)},
\end{aligned}
\]
where $d\m(t)=h(t)dt-d\n(t)$. Differentiating the modified
determinant $\p=D_2$  we have
\[
\lb{derfz} {\p'(k)\/\p(k)}=-2k\Tr \rt(R(k)-R_0(k)+R_0(k)VR_0(k)
\rt).
\]
Combining \er{derfx}, \er{derfz} we obtain \er{tre1}. If in addition
$V\in L^1(\R^3)$, then \er{B24} yields \er{tre2}. \BBox

We prove the second main result about the trace formulas.

\no {\bf  Proof of Theorem \ref{T3}.} We show \er{tr12}. Since
$\p\in \mH_\iy$  using  Theorem \ref{Tcf} and \er{Mas1} the
following asymptotic identity
$$
\begin{aligned}
\p(k)=\exp \rt[-i{Q_0+o(1)\/k } \rt] =\exp \rt[
-i{B_0+K_0+o(1)\/k}\rt] \exp \rt[i{C_M+o(1)\/k}\rt]
\end{aligned}
$$
as $\Im k\to \iy$, which yields \er{tr12}. We show \er{trj}. Let
$V\in W_{m+1}, m\ge 1$. From Lemma \ref{Th1} and from asymptotics
\er{apm}
we obtain
$$
\begin{aligned}
\p(k)=\exp \rt[-i{Q_0\/k }-i{Q_2\/k ^{3}}-i{Q_4\/k ^{5}}-\dots
  -i{Q_{2m}+o(1)\/k ^{2m+1}}  \rt]
\\
=\exp \rt[ -i{B_0\/k}-i{B_1\/2k^2}-i{B_2\/3k^3}-.....\rt] \exp
\rt[-i{K_0\/k}-i{K_1\/k^2}....\rt]
\\
\exp i\rt[{\cJ_0-iI_{0}\/k}+{\cJ_1-iI_{1}\/k^2}+...
+{\cJ_{2m}-iI_{2m}\/k^{2m+1}}+{M_{2m}(k)\/k^{2m+1}}\rt]
\end{aligned}
$$
which yields $ Q_j={B_j\/j+1}+K_j-\cJ_j+iI_j,\ j=0,1,... $. Thus we
have \er{trj}, since $I_j=\Im Q_j$.   \BBox

\

\no {\bf  Proof of Theorem \ref{T4}.}
Let $r>0$. We rewrite  \er{tr12} in the form:
\[
\begin{aligned}
\lb{tr5}  B_0+{\n(\R)\/\pi }=X_1+X_2, \qqq X_1={1\/\pi}\int_{-r}^r
\x(k)dk,\qq X_2={1\/\pi}\int_{\R\sm [-r,r]}\x(k)dk,
\end{aligned}
\]
where $\x(k)=\log |D_4(k+i0)|, \ k\in\R$. We estimate $X_2$. From
\er{DA2x}, \er{F2} we have
\[
\lb{p25}  \x(k)\le {17\/2}\|Y_0(k)\|_{\cB_4}^4\le
17{C_2^4\/2}{\|V\|^4\/|k|^{2}} \qqq \forall k\in \R,\qq
\]
which yields
\[
\lb{p26} X_2={1\/\pi}\int_{\R\sm [-r,r]}\!\! \!\!  \x(k)dk\le
{17\/2\pi}\int_{\R\sm [-r,r]}\|Y_0(k)\|_{\cB_4}^4dk \le
A\int_{k>r}{\|V\|^4\/|k|^{2}}dk={A\/ r}\|V\|^4.
\]
where $A=17{C_2^4\/\pi}$. We estimate $X_1$. Using the identity
\er{d2x},  we rewrite $X_1$ in the form
\[
\lb{d2X}
\begin{aligned}
X_1=X_{11}+X_{12}-X_{13},\qq \\
X_{11}={1\/\pi}\int_{-r}^r \log |\p(k)|dk,\qq
X_{1j}={1\/\pi}\int_{-r}^r \Re \p_j(k)dk,\qq
\end{aligned}
\]
where $\p_j={1\/j}\Tr Y_0^j(k), j=1, 2$. Estimate \er{D1} gives
\[
\lb{X1a} X_{11}={1\/\pi}\int_{-r}^r \log |\p(k)|dk\le
Br\|V\|_{3\/2}^2,\qqq B={2C_\bu\/\pi}.
\]
Estimates \er{epe} gives
\[
\lb{X1b} X_{12}={1\/\pi}\int_{-r}^r \Re \p_2(k)dk\le
{2r^{1\/2}\/\pi}\|\p_2\|\le
{r^{1\/2}\/(4\pi)^{4\/3}}\|V\|_{3\/2}\|V\|.
\]
Using the estimate $|\p_3(k)|\le
\|Y_0(k)\|_{\cB_2}^2\|Y_0(k)\|_{\cB_4}$ and  \er{F2} and \er{B22} we
obtain
\[
\lb{X1c} X_{13}\le
{1\/3\pi}\int_{-r}^r\|Y_0(k)\|_{\cB_2}^2\|Y_0(k)\|_{\cB_4}dk\le
2K\|V\|_{3\/2}^2 \|V\|\int_{0}^r{dk\/k^{1\/2}}=K\|V\|_{3\/2}^2
\|V\|r^{1\/2}
\]
where $K={8C_\bu C_2\/3\pi}$. Collecting \er{d2X}-\er{X1c} at
$\a=\|V\|, \b=\|V\|_{3\/2}$ we obtain for all $r>$:
\[
\lb{X12} X_1+X_2\le {A\a^4\/ r}+ B\b^2 r+ \a\b r^{1\/2}E, \qq
E=(4\pi)^{-{4\/3}}+K\b.
\]
The function ${A\a^4\/ r}+ B\b^2$ has minimum at
$r_o=(A/B)^{1\/2}{\a^2\/\b}$. Then this yields
\[
\lb{X12x} X_1+X_2\le 2(AB)^{1\/2}\a^2\b+ (A/B)^{1\/4}\a^2\b^{1\/2}E,
\qq E=(4\pi)^{-{4\/3}}+K\b.
\]

$$
X_1+X_2\le C\rt(\|V\|^4r^{-1}+\|V\|_{3\/2}^2r+ \|V\|_{3\/2}
\|V\|r^{1\/2} +\|V\|_{3\/2}^2 \|V\|r^{1\/2} \rt)
$$

at $r=\|V\|^{2\/3}$ we obtain
$$
X_1+X_2\le C\rt(\|V\|^4r^{-1}+\|V\|_{3\/2}^2r+ \|V\|_{3\/2}
\|V\|r^{1\/2} +\|V\|_{3\/2}^2 \|V\|r^{1\/2} \rt)
$$

which yields \er{eBs}. \BBox

\section {Appendix, analytic functions in the upper half-plane}
\setcounter{equation}{0}

We describe functions from the Hardy spaces. If $f\in \mH_\iy$, then

$\bu $ the function $B\in \mH_\iy$ with $\|B\|_{\mH_\iy}\le 1$ and
\[
\lb{B3} \lim_{v\to+0} B(u+iv)=B(u+i0), \qqq |B(u+i0)|=1 \qqq   \ a.
 e., u\in \R,
\]
\[
\lb{B4} \lim _{v\to 0}\int_\R\log |B(u+iv)|du=0.
\]

\begin{lemma}
\lb{TY2} Let $f\in \mH_\iy$ and let all its zeros $\{z_j\}$ in
$\C_+$ be uniformly bounded by $r_0$. Then

i) the coefficient $ B_n=2\sum_{j} \Im z_j^{n+1}, n\ge 0$ satisfies
\[
\lb{AB1} |B_n|\le 2\sum |\Im z_j^{n+1}|<\iy \qqq \qqq \forall \ n\ge
0,
\]
\[
\lb{AB2} |B_n|\le {\pi\/2}(n+1)r_0^{n}B_0 \qqq \qqq \forall \ n\ge
1.
\]
ii) The function $\log B(z)$ has an analytic continuation from
$\C_+\sm \{|z|<r_0\}$ into the domain $\{|z|>r_0\}$, where $r_0=\sup
|k_n|$ and has the following Taylor expansion
\[
\lb{AB3}
\log
B(z)=-{iB_0\/z}-{iB_1\/2z^2}-{iB_2\/3z^3}-....-{iB_{n-1}\/nz^n}-....
\]

\end{lemma}

\no {\bf  Proof.} i) Consider the function $F(w)=f(z(w))$, where the
conformal mapping $w:\C_+\to \dD$  is given by
$$
w=w(z)={i-z\/i+z},\qqq z\in \C_+,\qqq iw+zw=i-z,\qqq
z=z(w)=i{1-w\/1+w}.
$$
Thus $F\in \mH_\iy(\dD)$ and the zeros of $F$ have the form
$w_j=w(z_j)$. We have the identity
\[
1-|w|^2=1-{|i-z|^2\/|i+z|^2}=1-{(1-y)^2+x^2\/|i+z|^2}={4y\/|i+z|^2},
\]
which yields
\[
\sum (1-|w_j|^2)=\sum {4y_j\/|i+z_j|^2}.
\]
 Consider the function $f_1(\l)=f(\sqrt \l), \l\in \C_+$ and
$w(\l)$ is the conformal mapping. The function $f_1$ on $\C_+$ is
analytic and uniformly bounded, i.e., $f_1\in \mH_\iy(\C_+)$. The
function $f_1$ has the zeros $\l_j=z_j^2\in \C_+$ which satisfy
$$
\sum_{\Im \l_j>0} \Im \l_j =\sum_{\Im z_j^2>0}\Im z_j^2<\iy.
$$
If we apply similar arguments for the function $f(\sqrt \l), \l\in
\C_-$, we obtain
$$
\sum_{\Im \l_j<0} |\Im \l_j| =\sum_{\Im z_j^2<0}|\Im z_j^2|<\iy.
$$
Thus we have proved \er{AB1} for the case $n=2$. The proof for the case
$n\ge 3$ is similar.

 We have
$$
z_j=|z_j|e^{i\f_j}, \qqq \f_j=\ca \f_j^+ & {\rm if} \ \f_j\in (0,{\pi\/2}]\\
                                   \pi-\f_j^- & {\rm if} \  \f_j\in
                                   (\pi,{\pi\/2})\ac .
$$
This yields
\[
B_0=2\sum \Im z_j=2\sum |z_j|\sin \f_j< 2\sum_{j, \pm} |z_j|\f_j^\pm=A
\]
and
\[
B_0=2\sum \Im z_j=2\sum |z_j|\sin \f_j>{4\/\pi}\sum_{j, \pm}
|z_j|\f_j^\pm={2\/\pi}A.
\]
These estimates yield
\[
\sum \Im z_j^n=\sum |z_j|^n\sin n\f_j
\]
and then
\[
|\sum_{j} \Im z_j^n|\le \sum_{j} |z_j|^n |\sin n\f_j^\pm| \le
nr_0^{n-1}\sum_{j, \pm} |z_j|\f_j^\pm=nr_0^{n-1}A\le
{\pi\/2}nr_0^{n-1}B_0.
\]
which yields \er{AB2}.

ii) has an analytic continuation from $\C_+$ into the domain
$\{|k|>r_0\}$, where $r_0=\sup |k_n|$ and has the following Taylor
series

Let $c={z_j\/z}$ and $ \wt c={\ol z_j\/z}$.  We have
$$
b_j(z)=\log {z-z_j\/z-\ol z_j}=\log (1-c)-\log (1-\wt c)=-\sum_{n\ge
1}{2i\Im z_j^n\/nz^n }=b_{jm}(z)-\sum_{n\ge m}{2i\Im z_j^n\/nz^n },
$$
where $b_{jm}(z)=-\sum_{n=1}^{m-1}{2i\Im z_j^n\/nz^n }$ and
$$
|\log B(z)|\le |\sum_{j}\log {z-z_j\/z-\ol z_j}|\le
\sum_{j}\sum_{n\ge 1}{2|\Im z_j^n|\/n|z|^n },
$$
and
$$
\rt|\log B(z)-\sum_{j}b_{jm}(z)\rt| \le \sum_{j}\sum_{n\ge m}{2|\Im
z_j^n|\/n|z|^n } \le  \sum_{n\ge m}{r_0^{n-1}A\/|z|^n }=
{A\/r_0}\rt({r_0\/|z|}\rt)^m {1\/1-{r_0\/|z|}},
$$
which yields \er{AB3}. \BBox

We need the following results (see p 128 \cite{Ko98}):

$ \bu$ {\it If $u\in L^p(\R)$ and $p\in (1,\iy)$, then
\[
\lb{Ki1} F(k)={i\/\pi}\int_\R{u(t)dt\/k-t}\in \mH_p
\]
and $\Re F(k)\to u(t)$ as $ k\to t$ a.e in t (as $k\to t$
non-tangentially as usual).

$ \bu$ If $f\in \mH_p, p\ge 1$, then }
\[
\lb{Ki2} \log |f(t+i0)|\in L_{loc}^p(\R).
\]

\medskip

We need some results about functions from Hardy spaces. We begin
with asymptotics.

{\bf Definition.} {\it A function $h$ belongs to the class
$\gX_m=\gX_m(\R)$ if $h\in L_{real, loc}^1(\R)$ and has the form
\[
\lb{M1}
\begin{aligned}
 h(t)=P_m(t)+{h_m(t)\/t^{m+1}},\qqq
P_m(t)={I_{0}\/t}+{I_{1}\/t^2}+...+{I_{m}\/t^{m+1}},
\\
{h_m(\cdot)\/(1+|\cdot|)^a}\in L^1(\R),\qqq  h_m(t)=o(1)\qqq as \qqq
t\to \pm \iy
\end{aligned}
\]
for some real constants $I_0,I_1,....,I_m$ and integer $m\ge 0$ and
$a<1$.}

\

Note that if $h\in \gX_m $ for some $m\ge 0$, then there exist
finite integrals (the principal value):
\[
\begin{aligned}
\lb{M3}  \cJ_j={\rm v.o.}{1\/\pi}\lim_{s\to \iy}\int_{-s}^s
h_{j-1}(t)dt, \qqq h_{j}=t^{j+1}(h(t)-P_{j}(t)),\qq h_{-1}=h,
\end{aligned}
\]
for all $j=0,1,2,...,m-1$. For $h\in\gX_m $ we define the integrals
\[
\lb{dmM} M(k)={1\/\pi}\int_\R { h(t)\/k-t}dt,\qqq
M_m(k)={1\/\pi}\int_\R { h_m(t)\/k-t}dt,\qqq   k\in \C_+.
\]
In order to obtain trace formulas in Theorem \ref{T3} we need to
determine asymptotics of $M_m$.

\begin{lemma}
\lb{Th1} i) Let $h\in \gX_m$ for some $m\ge 0$. Then  the following
identity holds true:
\[
\begin{aligned}
\lb{M4} M(k)={1\/\pi}\int_\R {
h(t)\/k-t}dt={\cJ_0-iI_{0}\/k}+{\cJ_1-iI_{1}\/k^2}+...
+{\cJ_m-iI_{m}+M_m(k)\/k^{m+1}},
\end{aligned}
\]
for any $k\in \C_+$, where the real constants $\cJ_0,....,\cJ_{m}$
and the function $M_m$ are given by \er{dmM} and satisfies
\[
\lb{asMm} M_m(i\t)=o(1) \qqq \as \qq \t\to \iy.
\]

\end{lemma}


\no {\bf  Proof.} Consider the case  $m=0$, the proof for the case
$m\ge 1$ is similar. From \er{M1} we get
\[
\lb{h0} h\in L_{loc}^1(\R), \qqq h(t)={I_{0}+h_0(t)\/t},\qq
h_0(t)=o(1)\qqq \as \qqq t\to \pm \iy,
\]
for some real constant $I_0$. Introduce functions
$$
f(t)=h(t)-f_0(t),\qqq f_0(t)={I_{0} t\/t^2+1},\qqq  t\in \R,
$$
and note that $f(t)=h(t)-f_0(t)={o(1)\/t}$ as $ t\to \pm\iy$. Let
$k\in \C_+$. Then we have
\[
\lb{He1}
\begin{aligned}
 F_0(k)={1\/\pi}\int_\R   {f_0(t)\/k-t}dt
 ={I_{0}\/\pi}\int_\R {t\/(t^2+1)}{dt\/(k-t)}=-{iI_{0}\/k+i}.
\end{aligned}
\]
Thus the well-defined function $M(k)={1\/\pi}\int_\R {h(t)\/k-t}dt$
has the form
\[
\lb{Mz}
\begin{aligned}
 M(k)={1\/\pi}\int_\R {
h(t)\/k-t}dt=F_0(k)+F(k)=-{iI_{0}\/k+i}+F(k),
\\
 F(k)={1\/\pi}\int_\R {f(t)\/k-t}dt=
 {1\/\pi}\lim_{s\to \iy}\int_{-s}^s {f(t)\/k-t}dt
 ={\cJ_0\/k}+{ M_0(k)\/ k}+{I_{0}\/k(k+i)},
\end{aligned}
\]
since using ${1\/k-t}= {1\/k}+{t\/k(k-t)}$ and $h_0=th-I_0$,  we
have
\[
\lb{Mz2}
\begin{aligned} \int_{-s}^s {f(t)\/k-t}dt = {1\/k}\int_{-s}^s
f(t)dt+{1\/k}\int_{-s}^s {tf(t)\/k-t}dt\\
={1\/k}\int_{-s}^s h(t)dt+{1\/k}\int_{-s}^s
\rt(h_0(t)+{I_{0}\/(t^2+1)} \rt){dt\/k-t}
\end{aligned}
\]
and
\[
\lb{Mz3} {1\/\pi }\lim_{s\to \iy}\int_{-s}^s h(t)dt=\cJ_0,\qqq
{1\/\pi }\lim_{s\to \iy}\int_{-s}^s{h_0(t)\/k-t}dt=\int_\R
{h_0(t)\/k-t}dt=M_0(k),
\]
\[
{1\/\pi }\lim_{s\to \iy}\int_{-s}^s {I_{0}\/(t^2+1)}{dt\/(k-t)}=
{1\/\pi }\int_\R {I_{0}\/(t^2+1)}{dt\/(k-t)}={I_{0}\/k+i}.
\]
Collecting \er{He1}- \er{Mz} we obtain
\[
\begin{aligned}
\lb{Me1}
 M(k)=-{iI_{0}\/k+i}+F(k)=-{iI_{0}\/k+i}+{\cJ_0\/k}+
{ M_0(k)\/ k}+{I_{0}\/k(k+i)} = {\cJ_0-iI_{0}+M_0(k)\/k}.
 \end{aligned}
\]
where
$$
\begin{aligned}
M_0(k)={1\/\pi }\int_\R {h_0(t)\/k-t}dt,\qq h_0(t)=th(t)-I_0,\qq
\cJ_0={\rm v.p.}{1\/\pi}\int_\R h(t)dt.
\end{aligned}
$$
In order to show \er{asMm} we define a function
$g_k(t)={(1+|t|)^a\/k-t}, t\in\R$ and note that $\|g_k\|_\iy=o(1)$
as $k=i\t, \t\to \it$. Then we have
$$
\begin{aligned}
|M_0(k)|={1\/\pi }\rt|\int_\R {h_0(t)\/k-t}dt\rt|\le
{\|g_k\|_\iy\/\pi }\int_\R {|h_0(t)|\/(1+|t|)^a} dt=o(1),\qq  \as \
k=i\t,\  \t\to \iy,
\end{aligned}
$$
which yields \er{asMm}. \BBox

We describe the canonical factorization.

\begin{theorem}\label{Tcf}
Let a function $f\in\mH_p$ for some $p\ge 1$ and $f(k)=1+o(1)$ as
$|k|\to \iy$, uniformly with respect to ${\arg}\,k \in [0,\pi]$.
 Then $f$ has a canonical factorization in $\C_+$ given by
\[
\begin{aligned}
\lb{af1} f=f_{in}f_{out},\\
f_{in}(k)=B(k)e^{-iK(k)},\qqq  K(k)={1\/\pi}\int_\R {d\n(t)\/k-t}.
\end{aligned}
\]
$\bu$  $d\n(t)\ge 0$ is some singular compactly supported measure on
$\R$ and  satisfies
\[
\lb{sms}
\begin{aligned}
\n(\R)=\int_\R d\n(t)<\iy,\qqq
 \supp \n \ss [-r_c, r_c],
 \end{aligned}
\]
for some $r_c>0$. In particular, if $f$ is continuous in $\ol\C_+$.
Then
\[
\lb{ms}
 \supp \n \ss \{k\in \R: f(k)=0\}\ss [-r_c, r_c].
\]
$\bu$  The function $K(\cdot)$ has an analytic continuation from
$\C_+$ into the domain $\C\sm [-r_c, r_c]$   and has the following
Taylor series
\[
\lb{Kn} K(k)=\sum_{j=0}^\iy {K_j\/k^{j+1}},\qqq K_j={1\/\pi}\int_\R
t^jd\n(t).
\]
$\bu$ $B$ is the Blaschke product for $\Im k>0$ given by \er{Bk}.

\no $\bu$ Let in addition ${\log |f(\cdot)|\/(1+|\cdot|)^a}\in
L^1(\R)$ for some $a<1$. Then the outer factor  $f_{out}$ is  given
by
\[
\lb{Do2x} f_{out}(k)=e^{iM(k)},\qqq  M(k)= {1\/\pi}\int_\R {\log
|f(t)|\/k-t} dt,\qq k\in \C_+.
\]

\end{theorem}

{\bf Remark.} 1) These results are crucial to determine trace
formulas in Theorem \ref{T3}.

2) The integral $M(k), k\in \C_+$ in \er{Do2x} converges absolutely
since $f(t)=1+{O(1)\/t}$ as $t\to \pm\iy$.

{\bf Proof.} From \er{afm}  we deduce that each zero $\in \ol\C_+$
of $f\in \mH_p^0$  belongs  to the half-disc $\{k\in \ol \C_+:
|k|\le r_{c}\}$, for some $r_{c}>0$. It is well known (see p.119
\cite{Ko98}) that  the function $f(k), k\in \C_+$ has a standard
factorization $f=f_{in}f_{out}$, where  $f_{in}$ is the inner factor
 given by
\[
f_{in}(k)=e^{i\g+i\a k}B(k)e^{-iK(k)},\qqq  K(k)=
{1\/\pi}\int_\R\rt({1\/k-t}+{t\/t^2+1} \rt) d\n(t)
\]
with $\g\in \R, \a\ge 0$, and $d\n(t)\ge 0$ is a singular compactly
supported  measure on $\R$ such that $\int_\R{d\n(t)\/1+t^2}<\iy$
and $\supp \n \ss [-r_c,r_c]$. Here $B$ is the Blaschke product for
$\Im z>0$ given by \er{Bk}, since all zeros of $\p$ are uniformly
bounded. Without loss of generality (since $ \supp \n \ss  [-r_c,
r_c]$) we can write $K(k)$ in the form
\[
\lb{Di1x}   K(k)={1\/\pi}\int_{-r_c}^{r_c}{ d\n(t)\/k-t}.
\]
Due to \er{Di1x} we deduce that  the function $K$ is real on the set
$\R \sm [-r_c,r_c]$ . Thus the function $K$ has an analytic
extension from $\C_+$ into the whole cut plane $\C\sm [-r_c,r_c]$.
Moreover, $K$ has the Taylor series in the domain $\{|k|>r_c\}$. We
determine the Taylor series. For large $|k|$ we have
$$
\begin{aligned}
K(k)= {1\/k\pi}\int_{-r_c}^{r_c} {d\n(t)\/1-(t/k)} =\sum_{n\ge 0}
\int_\R {t^nd\n(t)\/k^{n+1}}
={K_0\/k}+{K_1\/k^2}+{K_2\/k^3}+{K_3\/k^4}+....
\qq \as \ |k|>r_c,\\
\end{aligned}
$$
where $K_n={1\/\pi}\int_{-r_c}^{r_c} {t^nd\n(t)}$, uniformly with
respect to ${\arg}\,k \in [0,2\pi]$, which yields \er{Kn}.

The function $f_{out}$ is the outer factor given by
\[
\lb{fout} f_{out}(k)=e^{iM(k)},\qqq  M(k)=
{1\/\pi}\int_\R\rt({1\/k-t}+{t\/t^2+1} \rt) \log |f(t)|dt.
\]
Consider $M(k)$. Due to \er{afm} we have $f(t)=1+O(1/t)$ as $t\to
\pm \iy$. Thus without loss of generality we can rewrite $M(k)$ in
the form \er{Do2x}.

It is well known (see p.119 \cite{Ko98}) that if $f$ is continuous
in $\ol\C_+$. Then \er{ms} holds true.\BBox

Let $f\in \mH_p$ for some $0<p\le \iy$. For integer $m\ge 0$ we say
$f$ belongs the class $ \mH_p^m=\mH_p^m(\C_+)$ if $f$ satisfies
\[
\lb{afm}
  \log f(k)={Q_0\/ik }+{Q_1\/ik ^{2}}+{Q_2\/ik ^{3}}+\dots
  +{Q_m+o(1)\/ik ^{m+1}},
\]
 as $|k |\to \iy$ uniformly in $\arg k \in [0,\pi]$, for some
constants  $Q_j\in \C$. We describe a canonical factorization of
functions from $\mH_p^m$.

\begin{theorem}\label{Tcftr}
Let a function $f\in \mH_p^m$ for some $m\ge 0, p\ge 1$ and let the
function $h(t)=\log |f(t)|, t\in \R$ belong to $\gX_m$. Then $f$ has
a canonical factorization $ f=f_{in}f_{out}$ in $\C_+$ given by
Theorem \ref{Tcf}, where  the function $M$ satisfies  for any $k\in
\C_+$ the following identity:
\[
\begin{aligned}
\lb{aM} M(k)={1\/\pi}\int_\R {
h(t)\/k-t}dt={\cJ_0-iI_{0}\/k}+{\cJ_1-iI_{1}\/k^2}+...
+{\cJ_m-iI_{m}+M_m(k)\/k^{m+1}},\\
\end{aligned}
\]
 where
\[
\begin{aligned}
\lb{aMx} M_m(k)={1\/\pi}\int_\R { h_m(t)\/k-t}dt,\qq
 \cJ_j=v.p.{1\/\pi}\int_\R h_{j-1}(t)dt,
\qqq h_{j-1}={I_j+h_j(t)\/t} ,\qqq
\end{aligned}
\]
\[
\lb{aM3} M_m(k)=o(1)\qq \as \qq \Im k\to \iy,
\]
$j=0,1,...,m,$ and  $h_{-1}=h$. Moreover, the following trace
formulas hold true:
\[
\lb{atr} B_j+K_j=\Re Q_j+\cJ_j, \qqq j=0,1,..,m.
\]
\end{theorem}

{\bf Proof.} From Lemma \ref{Th1} we deduce that relations
\er{aM}-\er{aM3} hold true. From Theorem \ref{Tcf}  we have
\[
\lb{apBKM} i\log f(k)=i\log B(k)+K(k)-M(k),
\]
where

 $\bu$ the  function $K$  has the following Taylor series \er{Kn},

 $\bu$ the
function $\log B(k) $ has the following Taylor series \er{B6},

$\bu$ the  function $M$  has asymptotics in  \er{aM}-\er{aM3},

$\bu$ the  function $f$  has  asymptotics given by \er{afm}.

Substituting all these asymptotics into the identity \er{apBKM} we
obtain
$$
\begin{aligned}
i\log f(k)={Q_0\/k}+{Q_1\/k^2}+{Q_2\/k^3}+...+{O(1)\/k^{m+2}}
\\=
\rt({B_0+K_0\/k}+{{B_1\/2}+K_1\/k^2}+{{B_2\/3}+K_2\/k^3}+... \rt)
-\rt({\cJ_0-iI_{0}\/k}+{\cJ_1-iI_{1}\/k^2}+{\cJ_2-iI_{2}\/k^3}+...\rt)
\end{aligned}
$$
 which gives \er{atr}.
\BBox

\footnotesize

\no {\bf Acknowledgments.} \footnotesize Evgeny Korotyaev is
grateful to Ari Laptev for stimulating discussions about the
Schr\"odinger operators with complex potentials.
 He is also grateful to Alexei Alexandrov (St. Petersburg) for
stimulating discussions and useful comments about Hardy spaces.


\footnotesize
\begin{thebibliography}
{999}\setlength{\itemsep}{-\parskip} \footnotesize


\bibitem[B66]{B66} V. S. Buslaev, The trace formulae and certain asymptotic
estimates of the kernel of the resolvent for the Schr{\"o}dinger
operator in three-dimensional space (Russian), Probl. Math. Phys.
No. I, Spectral Theory and Wave Processes, (1966) 82-101. Izdat.
Leningrad Univ. Leningrad.

\bibitem[C81]{C81}
Y. Colin de Verdi\`ere,  Une formule de traces pour l'op\'erateur de
Schr\"odinger dans ${\R}^3$, Ann. Sci. \'Ecole Norm. Sup. 14 (1981),
27--39.

\bibitem[DN02]{DN02} E. B. Davies; J. Nath, Schr\"odinger operators
 with slowly decaying potentials.  J. Comput. Appl. Math. 148(2002), 1--28.


\bibitem[DHK09]{DHK09} M. Demuth; M. Hansmann; G. Katriel,
On the discrete spectrum of non-selfadjoint operators.  J. Funct.
Anal. 257(2009), no. 9, 2742--2759.


\bibitem[DHK13]{DHK13} M. Demuth; M. Hansmann; G. Katriel,
Eigenvalues of non-selfadjoint operators: A comparison of two
approaches, in: Mathematical Physics, Spectral Theory and Stochastic
Analysis, Springer, 2013, 107--163.



\bibitem[F11]{F11}
R. L. Frank, Eigenvalue bounds for Schr\"odinger operators with
complex potentials, Bull. Lond. Math. Soc., 43 (2011), 745–-750.

\bibitem[F15]{F15} R. L. Frank, Eigenvalue bounds for Schr\"odinger
operators with complex potentials. III. Preprint (2015),
http://arxiv.org/pdf/1510.03411v1.pdf


 \bibitem[FLLS06]{FLLS06}
 R. L. Frank; A. Laptev; E. H. Lieb; R. Seiringer, Lieb–-Thirring inequalities for Schr\"odinger
operators with complex-valued potentials, Letters in Mathematical
Physics, 77 (2006),  309-–316.

\bibitem[FLS16]{FLS16} R. L. Frank; A. Laptev; O. Safronov, On the number of eigenvalues of
Schr\"odinger operators with complex potentials,  J. Lond. Math.
Soc. (2) 94 (2016), no. 2, 377–-390.




\bibitem[FS17]{FS17} R.L. Frank; J. Sabin,
Restriction theorems for orthonormal functions, Strichartz
inequalities, and uniform Sobolev estimates. Amer. J. of Math., to
appear

\bibitem    [G81] {G81}  J. Garnett, Bounded analytic functions,
Academic Press, New York, London, 1981.


\bibitem[G85]{G85}
L. Guillop\'e, Asymptotique de la phase de diffu\-sion pour
l'op\'erateur de Schr\"o\-dinger dans ${\R}^n$, S\'eminaire E.D.P.,
1984--1985, Exp. No. V, Ecole Polytechnique, 1985.


\bibitem[GK69] {GK69} I. Gohberg; M. Krein,
Introduction to the theory of linear nonselfadjoint operators.
Translated from the Russian, Translations of Mathematical
Monographs, Vol. 18 AMS, Providence, R.I. 1969.




\bibitem[IK12] {IK12}
H. Isozaki; E. Korotyaev, Inverse problems,
 trace formulae for discrete Schr\"odinger operators.
 Ann. Henri Poincaré 13 (2012), no. 4, 751–-788.

\bibitem[J90] {J90} A.  Jensen, High energy asymptotics for the total scattering phase
in potential scattering theory. Functional-analytic methods for
partial differential equations (Tokyo, 1989), 187–195, Lecture Notes
in Math., 1450, Springer, Berlin, 1990.

\bibitem[J92] {J92}
A. Jensen, High energy resolvent estimates for Schr\"odinger
operators. Ideas and methods in quantum and statistical physics
(Oslo, 1988), 254–260, Cambridge Univ. Press, Cambridge, 1992.



\bibitem[J06] {J06}
A. Jerbashian, Functions of A-Bounded Type in the Half-Plane. Vol.
4. Springer, 2006.


\bibitem[Ko98]{Ko98} Koosis, P. Introduction to $H_p$ spaces,
 Vol. 115 of Cambridge Tracts in Mathematic, 1998.

\bibitem   [K17]{K17} E. Korotyaev,
Trace formulas for Schr\"odinger operators on lattice,  preprint,
arXiv:1702.01388.


\bibitem   [KL16]{KL16} E. Korotyaev; A. Laptev,
Trace formulae for Schr\"odinger operators with complex-valued
potentials on cubic lattices, preprint,  arXiv:1609.09703.

\bibitem [KP04]{KP04}
E. Korotyaev; A. Pushnitski, A trace formula and high-energy
spectral asymptotics for   the perturbed Landau Hamiltonian. J.
Funct. Anal. 217 (2004), no. 1, 221–-248.

\bibitem [KP03]{KP03}   Korotyaev, E.; Pushnitski, A. Trace formulae
and high energy asymptotics for the Stark operator. Comm. Partial
Differential Equations 28 (2003), no. 3-4, 817--842.

\bibitem [KS17]{KS17} Korotyaev, E; Safronov, O.
Stark operators operators with complex potentials, preprint, 2017.



\bibitem[LS09] {LS09}
A. Laptev; O. Safronov, Eigenvalue estimates for Schr\"odinger
operators with complex potentials.  Comm. Math. Phys. 292 (2009),
no. 1, 29–-54.


\bibitem[LL97]{LL97}
E. Lieb; M. Loss, Analysis, AMS, Graduate Studies in Math., Vol. 14,
1997.

\bibitem[LT76]{LT76} E. Lieb, W. Thirring, Inequalities for the
moments of the eigenvalues of the Schr\"odinger Hamiltonian and
their relation to Sobolev inequalities. Studies in Mathematical
Physics. Princeton University Press, Princeton (1976), pp. 269--303.


\bibitem[P82]{P82}
G. Popov, Asymptotic behaviour of the scattering phase for the
Schr\"odinger operator, C. R. Acad. Bulgare Sci. 35 (1982), no. 7,
885--888.









\bibitem[R91]{R91}
D. Robert, Asymptotique \`a grande energie de la phase de diffusion
pour un potentiel, Asymptot. Anal. 3 (1991), 301--320.


\bibitem   [S10]{S10} O. Safronov, Estimates for eigenvalues of the
Schr\"odinger operator with a complex potential. Bull. Lond. Math.
Soc. 42 (2010), no. 3, 452--456.

 \bibitem   [Sa10]{Sa10}  O.  Safronov, On a sum rule for Schr\"odinger
 operators with complex potentials. Proc. Amer. Math. Soc. 138 (2010),
 no. 6, 2107--2112.




\end{thebibliography}
\end{document}